\documentclass[reqno]{amsart}
\usepackage{amsmath, amssymb, amstext, amscd, latexsym, dsfont, enumerate}
\usepackage[latin1]{inputenc}
\usepackage{mathrsfs}
\usepackage{textcomp}

\theoremstyle{plain}
\newtheorem{Theorem}{Theorem}
\newtheorem{Defi}{Definition}[section]
\newtheorem{Lemma}{Lemma}[section]

\numberwithin{equation}{section}

\theoremstyle{definition}
\newtheorem*{Rem}{Remark}

\newcommand{\calF}{\mathscr{F}}
\newcommand{\bbP}{\mathbb{P}}
\newcommand{\calH}{\mathcal{H}}
\newcommand{\calI}{\mathcal{I}}
\newcommand{\calD}{\mathcal{D}}
\newcommand{\sub}{\subset}
\newcommand{\dsR}{\mathds{R}}
\newcommand{\dsN}{\mathds{N}}
\newcommand{\bbE}{\mathbb{E}}
\newcommand{\sfH}{\textit{\textsf{H}}}
\newcommand{\dd}{\mathrm{d}}

\newcommand{\calB}{\mathcal{B}}
\newcommand{\calS}{\mathcal{S}}

\renewcommand{\Re}{\operatorname{Re}}

\renewcommand{\hat}{\widehat}
\renewcommand{\tilde}{\widetilde}

\newcommand{\oH}{\hspace*{0.39em}\raisebox{0.6ex}{\Large{\textdegree}}\hspace{-0.72em}\operatorname{H}}

\begin{document}
\title[Fractional White Noise Perturbations]{Fractional White Noise Perturbations of Parabolic Volterra Equations}

\author[S. Sperlich and M. Wilke]{Stefan Sperlich and Mathias Wilke}

\address{Martin-Luther-Universit\"{a}t Halle-Wittenberg, Institut f\"ur Mathematik, Theodor-Lieser-Str. 5, 06120 Halle, Germany}

%\email{jan.pruess@mathematik.uni-halle.de}
\email{stefan.sperlich@mathematik.uni-halle.de\quad \rm
(\emph{corresponding author})}
\email{mathias.wilke@mathematik.uni-halle.de}


%\date{June 26, 2008}


\maketitle

%\centerline{\emph{{Dedicated to the memory of Gunter Lumer.}}}
\vspace{0.5cm}

\begin{quote}\footnotesize
\textbf{Abstract.} Aim of this work is to extend the results of
Clément, Da Prato \& Prüss \cite{WNP} on the fractional white noise
perturbation with Hurst parameter $\sfH\in(0,1)$. We will obtain
similar results and it will turn out that the regularity of the
solution $u(t)$ increases with Hurst parameter
$\sfH$.\end{quote}\vspace{1cm}

\section{Introduction and Notations}
Stochastic differential equations play an important role in studying
random influences on deterministic systems. For this purpose usually
Wiener processes are considered. However, this is not adequate if
the chronological independence of the stochastic perturbations is
not sufficiently warranted. Therefore we make use of the concept of
a fractional Brownian motion, which was introduced by Mandelbrot \& van Ness
\cite{MaNe68}.\\

We are given a separable Hilbert space $\calH$ with norm
$|\cdot|_\calH$ and inner product $(\cdot|\cdot)_\calH$. Let $A$ be
a closed linear densely defined operator in $\calH$, and $b\in
L_1(\dsR_+)$ a scalar kernel. As in \cite{WNP} we consider the
integro-differential equation
\begin{equation}\label{prob1}
  \left\{\begin{split}\dot{u}(t)+ \int_0^t b(t-\tau)Au(\tau)\dd\tau &=Q^{1/2}\dot{B}^\sfH(t),\quad t\geq 0,\\
    u(0)&=u_0.\end{split}\right.
\end{equation}

Here the initial value $u_0$ is assumed to be an element of $\calH$. Moreover,
$B^\sfH$ is fractional Brownian motion in $\calH$ with Hurst
parameter $\sfH\in(0,1)$, with corresponding fractional white noise
$\dot{B}^\sfH$ and the operator $Q$ is of trace class (see Hypothesis \textbf{(B)} below).\\

Because problem (\ref{prob1}) is motivated from applications of
linear viscoelastic material behavior, we consider $G\sub\dsR^N$ to
be an open and bounded domain and the operator $-A$ to be an
elliptic differential operator like the Laplacian, the elasticity
operator, or the Stokes operator, together with appropriate boundary
conditions (e.g. Prüss \cite[Section I.5]{EvolInt}). In the
following we are particulary interested in
the case $\calH=L_2(G)$. \\

\noindent\textbf{Hypothesis (A).} {\itshape $A$ is an unbounded,
selfadjoint, positive definite operator in $\calH$ with compact
resolvent. Consequently, the eigenvalues $\mu_n>0$ of $A$ form a
nondecreasing sequence with $\lim_{n\to\infty}\mu_n=\infty$, the
corresponding
eigenvectors $(e_n)_{n\in\dsN}\sub\calH$ form an orthonormal basis of $\calH$.}\\

\noindent\textbf{Hypothesis (e).} {\itshape There is a constant
$C>0$ such that \[|e_n(\xi)|\leq C\quad\text{and}\quad|\nabla
e_n(\xi)|\leq
C\mu_n^{1/2},\] for all $n\in\dsN$ and all $\xi\in G$, where $\nabla$ denotes the gradient with respect to the variable $\xi$.}\\

\noindent\textbf{Hypothesis (b).} {\itshape $b\in L_1(\dsR_+)$ is
3-monotone, i.e. $b$ and $-\dot{b}$ are nonnegative, nonincreasing,
convex; in addition,
\begin{equation}\label{parabol}
  \lim_{t\to 0}\frac{\frac{1}{t}\int_0^t \tau
  b(\tau)\dd\tau}{\int_0^t -\tau\dot{b}(\tau)\dd\tau}<\infty.
\end{equation}}

Prüss proved in \cite[Section I.1]{EvolInt} that if \textbf{(A)} and
\textbf{(b)} are valid, the integrated version of problem
(\ref{prob1}) admits a resolvent $S(t)$ (which is strongly
continuous, uniformly bounded by 1, with
$\lim_{t\to\infty}|S(t)|_{\calB(\calH)}=0$ and $S\in
L_1(\dsR_+;\calB(\calH))$) such that the unique mild solution of
(\ref{prob1}) is given by the variation of parameters formula
\begin{equation}
  u(t)=S(t)u_0+\int_0^t S(t-\tau)f(\tau)\dd \tau,\quad t\geq 0,
\end{equation}
whenever $u_0\in\calH$ and $f\in L_{1,\text{loc}}(\dsR_+;\calH)$.\\

By means of the spectral decomposition of $A$, the resolvent family
$S(t)$ can be written explicitly as
\begin{equation}
  S(t)x=\sum_{n=1}^\infty s_n(t)(x|e_n)e_n,\quad t\geq 0,\quad
  x\in\calH,
\end{equation}
where the scalar functions $s_n(t)$ are the solutions of the scalar
problems
\begin{equation}\label{sk}
  \dot{s}_n(t)+\mu_n\int_0^t b(t-\tau)s_n(\tau)\dd\tau=0,\quad t\geq
  0,\quad s_n(0)=1.
\end{equation}

For the reader's convenience we repeat the definition of a scalar
fractional Brownian motion.

\begin{Defi}
    A real valued Gaussian process $\beta^\sfH:=\{\beta^{\sfH}(t)\}_{t\geq 0}$
    on a probability space
    $(\Omega,\mathscr{F},\mathbb{P})$ is called a fractional Brownian
    motion with Hurst parameter $\sfH\in(0,1)$ if for all $s,t\in\dsR_+$
    \begin{enumerate}[(i)]
      \item $\beta^\sfH(0)=0$,
      \item $\mathbb{E}\beta^{\sfH}(t)=0$,
      \item $\operatorname{Cov}[\beta^{\sfH}(t),\beta^{\sfH}(s)]=\frac{1}{2}
        \left(t^{2\sfH}+s^{2\sfH}-|t-s|^{2\sfH}\right)$.
    \end{enumerate}
\end{Defi}

Next we want to give an abstract formulation of the assumptions on
the covariance $Q$ and the fractional white noise $\dot{B}^\sfH$.\\

\noindent\textbf{Hypothesis (B).} {\itshape
$Q\in\mathscr{L}_1(\calH)$ is selfadjoint, positive semi-definite
and commutes with the operator $A$, i.e. there is a sequence
$(\gamma_n)\in\ell_1(\dsR_+)$, such that $Qe_n=\gamma_n e_n$ for all
$n\in\dsN$. $B^\sfH(t)$ is of the form
\begin{equation}\label{fBm}
  (B^\sfH(t)|x)=\sum_{n=0}^\infty\beta^\sfH_n(t)(x|e_n),\quad t\in\dsR, \quad
  x\in\calH,
\end{equation} where $\beta^\sfH_n$ are mutually
independent real valued fractional Brownian motions with Hurst
parameter $\sfH\in(0,1)$ on the probability space
$(\Omega,\calF,\bbP)$.}\\

Here, the symbol $\mathscr{L}_1(\calH)$ denotes the space of nuclear
operators on $\calH$. Note, that due to Hypothesis \textbf{(B)} the
operator $Q^{1/2}$ is well defined and belongs to
$\mathcal{L}_2(\calH)$; the space of Hilbert-Schmidt operators on
$\calH$. It is well known that $B^\sfH(t)$ (as in (\ref{fBm})) is
not a well defined $\calH$-valued random variable. However, due to
$B^\sfH(t):\Omega\to\calH_{Q^{-1/2}}$, where $\calH_{Q^{-1/2}}$ is
the completion of $\calH$ with respect to the norm
$|x|^2_{Q^{-1/2}}:=|Q^{-1/2}x|_\calH$, $x\in\calH$, the forcing
function $f$ is well defined since $Q^{1/2}B^\sfH(t)$ is a
mapping with values in $\calH$. This observation yields an alternative strategy how to avoid the appearance of the trace class operator $Q$, namely if in considering the space $\calH_{Q^{-1/2}}$ instead of $\calH$.\\

In the sequel an upper index $\langle t\rangle$, $t>0$ at a function
$f:\dsR_+\to\dsR$ means
\[f^{\langle t\rangle}(\tau):=\begin{cases}
  f(t-\tau) &: \tau\leq t;\\ 0 &: \tau>t.
\end{cases}\]
Moreover, we will make use of the theory of integration with respect
to fractional Brownian motions, which is provided by Pipiras and
Taqqu \cite{FBM-Int}. Hence we denote the fractional integral of
order $\alpha>0$ of a function $\phi$ by $\calI^\alpha \phi$,
precisely this means
\[(\calI^\alpha \phi)(r)=\frac{1}{\Gamma(\alpha)}\int_\dsR
\phi(\tau)(\tau-r)^{\alpha-1}_+\dd \tau,\quad r\in\dsR,\] where
$(x)_+:=\max\{0,x\}$. Recall that the Marchaud fractional derivative $\calD^{\alpha}$
is defined as the left inverse of $\calI^\alpha$ for $\alpha>0$ (see \cite[Page 111]{Samko}), i.e. for appropriate functions $\phi$ it holds that
\begin{equation}
  \calD^\alpha(\calI^\alpha \phi)\equiv\phi.
\end{equation}
Next we want to characterize the class of integrands $f$ with
respect to a fractional Brownian motion, such that the integral
$\int_\dsR f(\tau)\dd B^\sfH(\tau)$ is well defined. In order to
study the most general case, we consider the space $\Lambda_\sfH$
for integrands in the time domain which arises as
\begin{equation}
  \Lambda_\sfH:=\left\{f:\int_\dsR \left[(\calD^{\frac{1}{2}-\sfH}
f)(r)\right]^2\dd r<\infty\right\}\quad\text{for}\quad
0<\sfH<\tfrac{1}{2},
\end{equation}
or alternatively as
\begin{equation}
  \Lambda_\sfH:=\left\{f:\int_\dsR \left[(\calI^{\sfH-\frac{1}{2}}
f)(r)\right]^2\dd r<\infty\right\}\quad\text{for}\quad
\tfrac{1}{2}<\sfH<1.
\end{equation}
In both cases $\Lambda_\sfH$ is a linear space with inner
product\[(f|g)_{\Lambda_\sfH}=\frac{\Gamma^2(\sfH+\tfrac{1}{2})}{\zeta^2(\sfH-\tfrac{1}{2})}
\int_\dsR
(\calD^{\frac{1}{2}-\sfH}f)(r)(\calD^{\frac{1}{2}-\sfH}g)(r)\dd r,\]
or accordingly
\[(f|g)_{\Lambda_\sfH}=\frac{\Gamma^2(\sfH+\tfrac{1}{2})}{\zeta^2(\sfH-\tfrac{1}{2})}
\int_\dsR
(\calI^{\sfH-\frac{1}{2}}f)(r)(\calI^{\sfH-\frac{1}{2}}g)(r)\dd r,\]
where
\begin{equation*}
  \zeta(\sfH)=\left[\int_0^\infty
  \left[(1+\tau)^{\sfH}-\tau^{\sfH}\right]^2\dd\tau
  + \frac{1}{2\sfH+1}\right]^{1/2}.
\end{equation*}\\

\noindent\textbf{Remark.} In the literature (e.g.\ \cite{BHOZ08,Nua06}) the space $\Lambda_\sfH$ is often defined by a scalar product of the form
\[(f| g)_{\Lambda_\sfH}=c_\sfH\int_\dsR\int_\dsR f(s)g(t)|s-t|^{2\sfH-2}\dd s\dd t,\]
where $c_\sfH$ is an appropriate constant depending on $\sfH$. Recall that $f$ and $g$ might be distributions.\\

Pipiras and Taqqu proved in \cite[Proposition 3.2]{FBM-Int}, that
the embeddings
\begin{equation}\label{embeddings}
  L_1(\dsR)\cap L_2(\dsR)\hookrightarrow L_{1/\sfH}(\dsR)\hookrightarrow
  \Lambda_\sfH,
\end{equation}
hold true for $\sfH\in(\tfrac{1}{2},1)$.\\
In the spectral domain we are interested in integrands being a
member of the homogeneous Bessel potential space of order
$\tfrac{1}{2}-\sfH$,
\begin{equation*}
  \dot{\operatorname{H}}^{\frac{1}{2}-\sfH}_2(\dsR)=\left\{f\in \calS^*(\dsR):\int_\dsR
|\mathcal{F}f(\tau)|^2|\tau|^{-2\sfH+1}\dd\tau<\infty\right\},
\end{equation*}
where $\calS^*$ is the space of tempered distributions. It is well
known, that for
$f\in\dot{\operatorname{H}}^{\frac{1}{2}-\sfH}_2(\dsR)$ the Fourier
transform of $\calI^{\sfH-\frac{1}{2}}f$ or
$\calD^{\frac{1}{2}-\sfH}f$ is
\begin{equation}\label{Fourier}
  \psi_{\sfH-\frac{1}{2}}(x)(\mathcal{F}f)(x)|x|^{\frac{1}{2}-\sfH}
  =(\mathcal{F}f)(x)(ix)^{\frac{1}{2}-\sfH},
\end{equation}
where
\[\psi_{\alpha}(x)=e^{-i\pi \alpha/2}\chi_{\{x>0\}} +
e^{i\pi \alpha/2}\chi_{\{x<0\}},\quad x\in\dsR.\] Here $\chi_{M}$
denotes the indicator function of the set $M$. Hence by Plancherel's
Theorem it holds that
\begin{equation}\label{embedding2}
  \dot{\operatorname{H}}^{\frac{1}{2}-\sfH}_2(\dsR)\cong\Lambda_\sfH.
\end{equation}
An easy calculus shows that the identity
\begin{equation}\label{isometry}
  (f|g)_{\Lambda_\sfH}=\bbE\left[\left(\int_\dsR f(\tau)\dd
  \beta^\sfH(\tau)\right)\left(\int_\dsR g(\tau)\dd
  \beta^\sfH(\tau)\right)\right]
\end{equation}
holds for all $f,g\in\mathcal{E}$ , where $\mathcal{E}$ denotes the
set of all elementary functions. For a definition of the Wiener integral with respect to a fractional Brownian motion we refer to \cite{FBM-Int}. Since $\mathcal{E}$ is dense in
$\Lambda_\sfH$ (see \cite[Theorems 3.2 resp.\ 3.3]{FBM-Int}) equation
(\ref{isometry}) holds for all $f,g\in\Lambda_\sfH$.\\

\noindent\textbf{Remark.}
\begin{enumerate}[(i)]
  \item Since by (\ref{Fourier}) $\mathcal{F}(\calI^{-\kappa}f)\equiv
  \mathcal{F}(\calD^{\kappa}f)$ and by Plancherel's Theorem the norm in $\dot{\operatorname{H}}^\kappa_2(\dsR)$ can be
  rewritten as
  \begin{equation}\label{normh}
    |f|_{\dot{\operatorname{H}}^\kappa_2(\dsR)}=
    \begin{cases}
      \left|\calD^\kappa f\right|_{L_2(\dsR)} &: \kappa\geq 0;\\
      \left|\calI^{-\kappa} f\right|_{L_2(\dsR)} &: \kappa< 0.
    \end{cases}
  \end{equation}
  \item Observe that equation (\ref{isometry}) also holds on an arbitrary
  set $M\subset\dsR$. This can be seen by replacing
  $f$ and $g$ by $f\chi_M$ and $g\chi_M$, respectively.\\
\end{enumerate}

The plan of our paper is as follows. In Section 2 we state the main
results about fractional white noise perturbations of equations in
linear viscoelasticity, i.e. equation (\ref{prob1}), assuming the
Hypotheses \textbf{(A)}, \textbf{(b)}, and \textbf{(B)} explained
above. These results are proved in Section 3 by means of the methods
introduced in the monograph by Da Prato and Zabczyk
\cite{StochEqua}, adapted to evolutionary integral equations in
Clément and Da Prato \cite{CLE1}, \cite{CLE2}. The required
estimates were already available and taken from Monniaux and Prüss
\cite{Mon} and Clément, Da Prato \& Prüss \cite{WNP}.\\
Section 4 is devoted to a study of the equation \[u+g_\alpha*Au =
g_\beta*Q^{1/2}\dot{B}^\sfH\] on the halfline, where
$g_\kappa(t)=t^{\kappa-1}/\Gamma(\kappa)$, $t>0$ for $\kappa>0$
denotes the Riemann-Liouville kernel of fractional integration.\\

\section{Main results}
Concentrating on the stochastic case we let $h(t)=0$, i.e.
$f(t)=Q^{1/2}\dot{B}^\sfH(t)$; w.l.o.g. we set $u_0=0$. This means
that we have to investigate the stochastic convolution
\begin{equation}
  u(t)=\int_0^t S(t-\tau)\dd (Q^{1/2}B^\sfH)(\tau),\quad t\geq 0.
\end{equation}
In virtue of the spectral decompositions of $A$ and $Q$ we may
rewrite
\begin{equation}\label{sprob1}
  u(t)=\sum_{n=1}^\infty \sqrt{\gamma_n}\int_0^t s_n(t-\tau)e_n\dd
  \beta^\sfH_n(\tau),\quad t\geq 0.
\end{equation}
Our main result on problem (\ref{prob1}) reads as follows.\\

\begin{Theorem}\label{Theorem}
Let $\sfH\in(0,1)$. Assume that Hypotheses \textbf{(A)},
\textbf{(b)}, \textbf{(B)} are valid and suppose
\begin{equation}
  \sum_{n=1}^\infty
  \gamma_n \mu_n^{-\frac{2\sfH}{\rho}}<\infty,
\end{equation}
where
\begin{equation}\label{rho}
  \rho:=1+\frac{2}{\pi}\sup\{|\operatorname{arg}\hat{b}(\lambda)|:\Re\lambda>0\}.
\end{equation}
Then the series (\ref{sprob1}) converges in $L_2(\Omega;\calH)$,
uniformly in $t$ on bounded subsets of $\dsR_+$ and $u\in
C_b(\dsR_+;L_2(\Omega;\calH))$. $u(t)$ is a Gaussian random variable
with mean zero and covariance operator $Q_t$, defined by
\begin{equation}
  Q_tx=\sum_{n=1}^\infty \left\|s_n^{\langle
  t\rangle}\sqrt{\gamma_n}\right\|^2_{\Lambda_\sfH}(x|e_n)
  e_n,\quad x\in\calH,
\end{equation}
and we have $\operatorname{Tr}[Q_t]\leq
c_\sfH\operatorname{Tr}[QA^{-2\sfH/\rho}]$.\\
If in addition, there is
$\theta\in(0,1)$ such that
\begin{equation}
  \sum_{n=1}^\infty
  \gamma_n\mu_n^{\frac{2\sfH(\theta-1)}{\rho}}<\infty,
\end{equation}
then $u\in C_b^{\alpha}(\dsR_+;L_2(\Omega;\calH))$ for each $\alpha\in(0,\theta \sfH)$.\\
In case $\calH=L_2(G)$ and Hypothesis \textbf{(e)} as well as
\begin{equation}
  \sum_{n=1}^\infty
  \gamma_n \mu_n^{\theta-\frac{2\sfH}{\rho}}<\infty,
\end{equation}
are met, then $u\in C_b(\dsR_+;C^\alpha(G;L_2(\Omega)))$ for each $\alpha\in(0,\theta)$.
\end{Theorem}

Here and in the sequel we denote by $c_\sfH>0$ a generic constant
depending on $\sfH$.\\

\noindent\textbf{Remark.} \begin{enumerate}[(i)]
\item $\Lambda_{\frac{1}{2}}$ is isometrically isomorphic to $L_2(\dsR)$. In
this sense Theorem \ref{Theorem} is a generalization of
\cite[Theorem 2.1]{WNP}.
\item Note that by Hypothesis \textbf{(A)} and \textbf{(b)} the problem
under consideration is \emph{parabolic}, i.e. $\rho\in[1,2)$.\\
\end{enumerate}

\section{Proof of the main results}
The idea of the proof is, of course, similar to Clément et al.
\cite{WNP} and follows the arguments for the Cauchy problem
presented
in Da Prato and Zabczyk \cite{StochEqua}.\\

Let us cite a useful lemma which was proven in \cite[Lemma
3.1]{WNP}:

\begin{Lemma}\label{Lemma_sn}
  Suppose the kernel $b(t)$ is subject to Hypothesis \textbf{(b)},
  and let $\rho\in(1,2)$ be defined by (\ref{rho}). Then for every
  $n\in\dsN$ it is
  \begin{enumerate}[(i)]
    \item $|s_n(t)|\leq 1$ for all $t,\mu_n>0$;
    \item $|\dot{s}_n|_{L_1(\dsR_+)}\leq C$ for all $\mu_n>0$;
    \item $|t\dot{s}_n|_{L_1(\dsR_+)}\leq C\mu_n^{-1/\rho}$ for all $\mu_n>0$;
    \item $|s_n|_{L_1(\dsR_+)}\leq C\mu_n^{-1/\rho}$ for all $\mu_n>0$,
  \end{enumerate}
  where $C>0$ denotes a constant which is independent of $\mu_n>0$.
\end{Lemma}

Now, let the hypotheses of Theorem \ref{Theorem} be fulfilled.
Observe that for $\sfH\in(\tfrac{1}{2},1)$ by (iv) and (i) of Lemma
\ref{Lemma_sn} the functions $s^{\langle t\rangle}_n$ belongs to
$L_1(\dsR)\cap L_2(\dsR)$ and hence by embedding (\ref{embeddings})
to $\Lambda_\sfH$. So by identity (\ref{isometry}) we obtain
\begin{equation}\label{Eu}
  \bbE|u(t)|_\calH^2=c_\sfH\sum_{n=1}^\infty
  \gamma_n \int_\dsR\left[(\calI^{\sfH-\frac{1}{2}} s_n^{\langle t\rangle})(r)\right]^2\dd r=c_\sfH\sum_{n=1}^\infty
  \gamma_n\left|\calI^{\sfH-\frac{1}{2}} s_n^{\langle t\rangle}\right|^2_{L_2(\dsR)}.
\end{equation}
As a result of \cite[Theorem 5.3]{Samko} the operator
$\calI^{\sfH-\frac{1}{2}}$ is bounded from $L_{1/\sfH}(\dsR)$ into
$L_2(\dsR)$. Thus we have by (i) and (iv) of Lemma \ref{Lemma_sn}
\begin{equation}\label{Eu2}
  \bbE|u(t)|_\calH^2\leq c_\sfH\sum_{n=1}^\infty\gamma_n
  |s_n^{\langle t\rangle}|^2_{L_{1/\sfH}(\dsR)}\leq c_\sfH\sum_{n=1}^\infty\gamma_n
  |s_n^{\langle t\rangle}|^{2\sfH}_{L_1(\dsR)} \leq
  c_\sfH\sum_{n=1}^\infty\gamma_n
  \mu_n^{-2\sfH/\rho}
\end{equation}
which is finite by assumption. In the case $\sfH\in(0,\tfrac{1}{2})$
one may argue as in the latter situation to obtain with the aid of (\ref{normh}) and the continuous embedding $\operatorname{H}^\kappa_2(\dsR)\hookrightarrow \dot{\operatorname{H}}^\kappa_2(\dsR)$, $\kappa>0$,
\begin{equation}
  \bbE|u(t)|_\calH^2=c_\sfH\sum_{n=1}^\infty
  \gamma_n\left|\calD^{\frac{1}{2}-\sfH} s_n^{\langle
  t\rangle}\right|^2_{L_2(\dsR)}\leq c_\sfH\sum_{n=1}^\infty
  \gamma_n|s_n^{\langle
  t\rangle}|^2_{\operatorname{H}^{\frac{1}{2}-\sfH}_2(\dsR)}.
\end{equation}
Hence by interpolation and Lemma \ref{Lemma_sn}
\begin{equation}
  \bbE|u(t)|^2_\calH\leq c_\sfH\sum_{n=1}^\infty
  \gamma_n|s_n^{\langle
  t\rangle}|^{2\sfH}_{L_1(\dsR)}\cdot|s_n^{\langle
  t\rangle}|^{2(1-\sfH)}_{\operatorname{H}^1_1(\dsR)}\leq c_\sfH\sum_{n=1}^\infty\gamma_n
  \mu_n^{-2\sfH/\rho}
\end{equation}
holds. This can be seen as follows. Let us denote by $[X;Y]_\delta$ the
complex interpolation space of the spaces $X$ and $Y$ with
parameter $\delta\in(0,1)$. Then $[X;Y]_\delta=Z$ entails the interpolation
inequality $|f|_Z\leq c|f|_X^{1-\delta}|f|_Y^{\delta}$ for all
$f\in Z$. It follows from \cite[Theorem 2.4.7]{Tri83} that
\begin{equation*}
  \left[L_1(\dsR_+);\operatorname{H}^\tau_{q}(\dsR_+)\right]_\delta=\operatorname{H}^{\frac{1}{2}-\sfH}_2(\dsR_+)
\end{equation*}
holds with
\begin{equation}\label{interpol}
  \delta=1-\sfH,\quad \tau=\frac{1-2\sfH}{2(1-\sfH)},\quad
  q=\frac{2(1-\sfH)}{1-2\sfH}.
\end{equation}
This observation yields
\begin{equation*}
  |s_n^{\langle
  t\rangle}|_{\dot{\operatorname{H}}^{\frac{1}{2}-\sfH}_2(\dsR_+)}\leq |s_n^{\langle
  t\rangle}|_{\operatorname{H}^{\frac{1}{2}-\sfH}_2(\dsR_+)}
  \leq
  c |s_n^{\langle
  t\rangle}|^\sfH_{L_1(\dsR_+)}|s_n^{\langle
  t\rangle}|^{1-\sfH}_{\operatorname{H}^\tau_q(\dsR_+)}.
\end{equation*}
Now, one may apply \cite[Theorem 2.7.1]{Tri83} to verify that the
embedding
\begin{equation*}
  \operatorname{H}^1_1(\dsR_+) \hookrightarrow
  \operatorname{H}^\tau_q(\dsR_+)
\end{equation*}
holds with $\tau$ and $q$ as in (\ref{interpol}).\\

Thus $u(t)$ is a zero mean $\calH$-valued Gaussian random variable.
Let $Q_t$ be its covariance operator, then for
$\sfH\in(\tfrac{1}{2},1)$
\begin{equation*}\begin{split}
  (Q_t x|y)_\calH&=\bbE\left[(u(t)|x)(u(t)|y)\right]\\
  &= \sum_{n=1}^\infty \gamma_n(e_n|x)(e_n|y) \bbE\left|\int_\dsR
  s_n^t(\tau)\dd\beta^\sfH_n(\tau)\right|^2\\
  &= \frac{\Gamma^2(\sfH+\tfrac{1}{2})}{\zeta^2(\sfH-\tfrac{1}{2})}\sum_{n=1}^\infty
  (e_n|x)(e_n|y) \int_\dsR \left[(\calI^{\sfH-\frac{1}{2}}
  s_n^{\langle t\rangle}\sqrt{\gamma_n})(\tau)\right]^2\dd\tau\\
  &=\frac{\Gamma^2(\sfH+\tfrac{1}{2})}{\zeta^2(\sfH-\tfrac{1}{2})}
  \sum_{n=1}^\infty\left(\int_\dsR
  \left[(\calI^{\sfH-\frac{1}{2}}s_n^{\langle
  t\rangle}\sqrt{\gamma_n})\right]^2(\tau)\dd\tau(e_n|x)e_n\Big|y\right),
\end{split}\end{equation*}
and with help of (\ref{Eu}) and (\ref{Eu2})
$\operatorname{Tr}[Q_t]\leq
c_\sfH\operatorname{Tr}[QA^{-2\sfH/\rho}]$ follows. Replacing
$\calI^{\sfH-\frac{1}{2}}$ by $\calD^{\frac{1}{2}-\sfH}$ yields the
claim for
$\sfH\in(0,\tfrac{1}{2})$.\\

Concerning Hölder-continuity we will use the following two estimates
with the convention $s_n(\tau)=0$ for $\tau<0$.

\begin{Lemma}\label{Lemma_est}
  Suppose that the kernel $b(t)$ is subject to Hypothesis
  \textbf{(b)} and let $\kappa\in(1,2)$. Then for each $\theta\in(0,1)$ there is a constant
  $C_\theta>0$ such that
  \begin{equation}
    \int_x^t |s_n(t-\tau)|^{\kappa}\dd \tau\leq
    C_\theta\mu_n^{(\theta-1)/\rho}|t-x|^\theta,\quad 0<x<t,
  \end{equation}
  and
  \begin{equation}
    \int_{-\infty}^x |s_n(t-\tau)-s_n(x-\tau)|^{\kappa}\dd\tau
    \leq C_\theta\mu_n^{(\theta-1)/\rho}|t-x|^\theta,\quad x<t.
  \end{equation}
\end{Lemma}
The proof of Lemma \ref{Lemma_est} follows exactly the lines of
\cite[Proof of Lemma 3.1]{WNP}. Therefore we omit it. For
$\sfH\in(\tfrac{1}{2},1)$ we use the identity (\ref{isometry}) to
obtain
\[\begin{split}
  \bbE|u(t)-u(x)|^2_\calH&=\bbE(u(t)-u(x)|u(t)-u(x))_\calH =
  \sum_{n=1}^\infty \gamma_n (s_n^{\langle t\rangle}-s_n^{\langle x\rangle}|s_n^{\langle t\rangle}-s_n^{\langle x\rangle})_{\Lambda_\sfH}\\
  &= c_\sfH \sum_{n=1}^\infty \gamma_n
  \left|\calI^{\sfH-\frac{1}{2}}(s_n^{\langle t\rangle}-s_n^{\langle x\rangle})\right|^2_{L_2(\dsR)}\leq c_\sfH \sum_{n=1}^\infty
  \gamma_n |s_n^{\langle t\rangle}-s_n^{\langle x\rangle}|^2_{L_{1/\sfH}(\dsR)}
\end{split}\]
and we have
\[\begin{split}|s_n^{\langle t\rangle}-s_n^{\langle x\rangle}|_{L_{1/\sfH}}&= \left[\int_{-\infty}^x
|s_n(t-\tau)-s_n(x-\tau)|^{1/\sfH}\dd\tau
+\int_x^t|s_n(t-\tau)|^{1/\sfH}\dd\tau\right]^\sfH.
\end{split}\] For $\sfH\in(0,\tfrac{1}{2})$ it is
\begin{equation*}
  \bbE|u(t)-u(x)|^2_\calH=c_\sfH \sum_{n=1}^\infty \gamma_n
  \left|\calD^{\frac{1}{2}-\sfH}(s_n^{\langle t\rangle}-s_n^{\langle x\rangle})\right|^2_{L_2(\dsR)}
\end{equation*}
and the estimate
\begin{equation}
  \left|\calD^{\frac{1}{2}-\sfH}(s_n^{\langle t\rangle}-s_n^{\langle
  x\rangle})\right|^2_{L_2(\dsR)} \leq c |s_n^{\langle t\rangle}-s_n^{\langle
  x\rangle}|^2_{\operatorname{H}^{\frac{1}{2}-\sfH}_2(\dsR)}\leq \tilde{c} |s_n^{\langle t\rangle}-s_n^{\langle
  x\rangle}|^{2\sfH}_{L_1(\dsR)}
\end{equation}
holds for sufficient large $n\in\dsN$, by interpolation and Lemma
\ref{Lemma_sn}. Thus by employing Lemmata \ref{Lemma_sn} and
\ref{Lemma_est} this yields
\[\bbE|u(t)-u(x)|^2_\calH\leq c_\sfH|t-x|^{2\theta \sfH}\sum_{n=1}^\infty \gamma_n
\mu_{n}^{2\sfH(\theta-1)/\rho},\] for $\sfH\in(0,1)$ and with the
aid of Kahane-Khinchine inequality (e.g. \cite[Corollary
3.4.1]{KwWo92}) and the Kolmogorov-\v{C}entsov-Theorem (e.g.
\cite[Theorem 2.8]{KaSh91}) we may conclude Hölder-continuity with
respect to $t$ of $u(t)$ as in the proofs given in Clément and Da
Prato \cite{CLE1} or \cite{CLE2}. Similarly, in case \textbf{(e)}
holds, we obtain spatial Hölder-continuity from the identities
\[\bbE|u(t,\xi)-u(t,\eta)|^2=c_\sfH\sum_{n=1}^\infty
\gamma_n \left|\calI^{\sfH-\frac{1}{2}} s_n^{\langle
t\rangle}\right|^2_{L_2(\dsR)} |e_n(\xi)-e_n(\eta)|^2\] for
$\sfH\in(\tfrac{1}{2},1)$ and
\[\bbE|u(t,\xi)-u(t,\eta)|^2=c_\sfH\sum_{n=1}^\infty
\gamma_n \left|\calD^{\frac{1}{2}-\sfH} s_n^{\langle
t\rangle}\right|^2_{L_2(\dsR)} |e_n(\xi)-e_n(\eta)|^2\] for
$\sfH\in(0,\tfrac{1}{2})$ respectively.\\

\section{Fractional derivatives and fractional white noise}
In the remaining part of this paper we take up a different viewpoint
to equations with fractional noise. We consider the problems
\begin{equation}\label{fracprob}
  u+g_\alpha*Au=g_\beta*Q^{1/2}\dot{B}^\sfH
\end{equation}
in the Hilbert space $\calH$, where the operator $A$ is subject to
Hypothesis \textbf{(A)} and also to \textbf{(e)} if appropriate, the
covariance $Q$ and the fractional Brownian motion $B^\sfH$ are
subject to \textbf{(B)}, and $g_\kappa$ denotes the fractional
integration kernel
\[g_\kappa(t)=\frac{t^{\kappa-1}}{\Gamma(\kappa)},\quad t>0,\]
where $\kappa>0$. Note that the kernel $g_\alpha$ is of
subexponential growth, i.e.
\begin{equation*}
  \int_0^\infty e^{-\omega t}|g_\alpha(t)|\dd t<\infty
\end{equation*}
for arbitrary small $\omega>0$. This means that that Laplace
transform $\hat{g}_\alpha$ is well defined.

\noindent\begin{Rem}
  Problem (\ref{fracprob}) with $\beta=1$, has also been studied in
  a recent paper of Bonaccorsi \cite{Bon07} in regard to existence of a mild solution.
\end{Rem}

For $\alpha\in(0,2)$, $\beta>0$, define the scalar fundamental
solution of (\ref{fracprob}) by
\begin{equation}
  \hat{r}_n(\lambda)=\frac{\hat{g}_\beta(\lambda)}{1+\mu_n\hat{g}_\alpha(\lambda)}=\frac{\lambda^\alpha}{\lambda^\beta(\lambda^\alpha
  +\mu_n)},\quad \Re\lambda>0,\quad \mu_n>0,
\end{equation}
where $\hat{r}_n$ denotes the Laplace transform of $r_n$.
Furthermore with the convention $r_n(\tau)=0$ for $\tau<0$ we have
by the Paley-Wiener Theorem
\begin{equation}\label{rnnorm}\begin{split}
  |r_n|^2_{\dot{\operatorname{H}}^{\frac{1}{2}-\sfH}_2(\dsR)}&= \int_\dsR
  |(\mathcal{F}r_n)(\rho)|^2|\rho|^{1-2\sfH}\dd \rho\\
  &\leq c_\alpha\int_\dsR
  \Bigl[\frac{|\rho|^\alpha}{|\rho|^\beta(|\rho|^\alpha+\mu_n)}\Bigr]^2|\rho|^{1-2\sfH}\dd\rho\\
  &=2c_\alpha \int_0^\infty
  \Bigl[\frac{\rho^\alpha}{\rho^\beta(\rho^\alpha+\mu_n)}\Bigr]^2\rho^{1-2\sfH}\dd\rho\\
  &= 2c_\alpha \mu_n^{\frac{2(1-\beta-\sfH)}{\alpha}}\int_0^\infty
  \Bigl[\frac{\tau^{\alpha-\beta-\sfH+\frac{1}{2}}}{1+\tau^\alpha}\Bigr]^2\dd\tau,
\end{split}\end{equation} and the right integral is finite if and only if
$1-\sfH<\beta<1-\sfH+\alpha$. Thus by isomorphism (\ref{embedding2})
$r_n$ belongs to $\Lambda_\sfH$ whenever $\alpha\in(0,2)$ and
$\beta\in(1-\sfH,1-\sfH+\alpha)$. The solution of (\ref{fracprob})
can be rewritten as
\begin{equation}
  u(t)=\sum_{n=1}^\infty \sqrt{\gamma_n}\int_0^t r_n(t-\tau)\dd
  \beta^\sfH_n(\tau)e_n,\quad t>0,
\end{equation}
and therefore as in Section 3 it is by means of representation
(\ref{normh})
\begin{equation}\label{id1}
  \bbE|u(t)|_\calH^2=c_\sfH\sum_{n=1}^\infty
  \gamma_n|r_n^{\langle
  t\rangle}|^2_{\dot{\operatorname{H}}^{\frac{1}{2}-\sfH}_2(\dsR)}
\end{equation}
as well as
\begin{equation}\label{id2}
  \bbE|u(t)-u(x)|^2_\calH= c_\sfH\sum_{n=1}^\infty \gamma_n|r_n^{\langle
  t\rangle}-r_n^{\langle x\rangle}|^2_{\dot{\operatorname{H}}^{\frac{1}{2}-\sfH}_2(\dsR)}
\end{equation}
and in case $\calH=L_2(G)$ and \textbf{(e)} is valid
\begin{equation}
  \bbE|u(t,\xi)-u(t,\eta)|^2=c_\sfH\sum_{n=1}^\infty
\gamma_n |r_n^{\langle
t\rangle}|^2_{\dot{\operatorname{H}}^{\frac{1}{2}-\sfH}_2(\dsR)}
|e_n(\xi)-e_n(\eta)|^2.
\end{equation}
Moreover, it is due to
\begin{multline*}
  (\mathcal{F}f^{\langle t\rangle})(\xi)=\int_\dsR f(t-\tau)\chi_{(-\infty,t]}(\tau)e^{-i\xi\tau}\dd\tau\\
  =\int_\dsR f(-s)\chi_{(-\infty,0]}(s)e^{-i\xi(s+t)}\dd s=e^{-i\xi t}(\mathcal{F}f^{\langle 0\rangle})(\xi)
\end{multline*}
that for all $t\in\dsR$
\begin{equation}\label{H_shift}
  \|f^{\langle
  t\rangle}\|_{\dot{\operatorname{H}}^{\sigma}_{2}(\dsR)}=\|f^{\langle
  0\rangle}\|_{\dot{\operatorname{H}}^{\sigma}_{2}(\dsR)}=\|f\|_{\dot{\operatorname{H}}^{\sigma}_{2}(\dsR_+)}.
\end{equation}
Thus
 \[|r_n^{\langle
t\rangle}|_{\dot{\operatorname{H}}^{\frac{1}{2}-\sfH}_2(\dsR_+)}\leq|r_n^{\langle
t\rangle}|_{\dot{\operatorname{H}}^{\frac{1}{2}-\sfH}_2(\dsR)}=|r_n^{\langle
0\rangle}|_{\dot{\operatorname{H}}^{\frac{1}{2}-\sfH}_2(\dsR)}=|r_n|_{\dot{\operatorname{H}}^{\frac{1}{2}-\sfH}_2(\dsR_+)}\]
holds for $t\geq 0$, as soon as
$r_n\in\dot{\operatorname{H}}^{\frac{1}{2}-\sfH}_2(\dsR_+)$.
Identities (\ref{id1}) and (\ref{id2}) show that the solution $u(t)$
of (\ref{fracprob}) exists and is continuous in $L_2(\Omega;\calH)$
if and only if
\begin{equation}\label{sigma1}
  \sigma_1:=\sum_{n=1}^\infty \gamma_n|r_n|^2_{\dot{\operatorname{H}}^{\frac{1}{2}-\sfH}_2(\dsR_+)}<\infty.
\end{equation}
Next observe that we have for $\sfH\in(\tfrac{1}{2},1)$
\[\begin{split}\left|\calI^{\sfH-\frac{1}{2}}(r_n^{\langle t\rangle}-r_n^{\langle
x\rangle})\right|^2_{L_2(\dsR)}&= \int_\dsR
\left|(\calI^{\sfH-\frac{1}{2}} r_n^{\langle
t\rangle})(\tau)-(\calI^{\sfH-\frac{1}{2}} r_n^{\langle
x\rangle})(\tau)\right|^2\dd \tau\\
&=\int_\dsR \left|(\calI^{\sfH-\frac{1}{2}} r_n^{\langle
0\rangle})(\tau-t)-(\calI^{\sfH-\frac{1}{2}} r_n^{\langle
0\rangle})(\tau-x)\right|^2\dd \tau\\
&=\left|(\calI^{\sfH-\frac{1}{2}}r_n^{\langle 0\rangle})(x-t+\cdot)
- (\calI^{\sfH-\frac{1}{2}}r_n^{\langle
0\rangle})(\cdot)\right|^2_{L_2(\dsR)}\\
&\leq \left|\calI^{\sfH-\frac{1}{2}}r_n^{\langle
0\rangle}\right|^2_{B^\theta_{2,\infty}(\dsR)}|t-x|^{2\theta}
\end{split}\]
and analogue for $\sfH\in(0,\tfrac{1}{2})$
\begin{equation*}
  \left|\calD^{\frac{1}{2}-\sfH}(r_n^{\langle t\rangle}-r_n^{\langle
  x\rangle})\right|^2_{L_2(\dsR)}\leq
  \left|\calD^{\frac{1}{2}-\sfH}r_n^{\langle
  0\rangle}\right|^2_{B^\theta_{2,\infty}(\dsR)}|t-x|^{2\theta},
\end{equation*}
where $B^\theta_{2,\infty}(\dsR)$ denotes a Besov space, with the equivalent norm
\[|f|_{B^\theta_{2,\infty}(\dsR)}=\left[|f|^2_{L_2(\dsR)}+\sup_{h\in\dsR}\int_\dsR\frac{|f(y+h)-f(y)|^2}{|h|^{2\theta}}\dd y\right]^{1/2}.\]
Now we have
the embedding \[\operatorname{H}^\theta_{2}(\dsR) \hookrightarrow
B^\theta_{2,\infty}(\dsR),\] cf. \cite[Theorem 2.3.2 (c)]{Triebel},
and the apparent relation
\begin{equation}
  |f|_{\dot{\operatorname{H}}^\kappa_2(\dsR)} +
  |f|_{\dot{\operatorname{H}}^{\theta+\kappa}_2(\dsR)}=
  \begin{cases}
    \left|\calD^\kappa f\right|_{\operatorname{H}^\theta_2} &: \kappa\geq 0,\\
    \left|\calI^{-\kappa} f\right|_{\operatorname{H}^\theta_2} &: \kappa< 0.\\
  \end{cases}
\end{equation}
So the condition
\begin{equation}\label{sigma2}
  \sigma_2:=\sum_{n=1}^\infty\gamma_n\left[|r_n|_{\dot{\operatorname{H}}^{\frac{1}{2}-\sfH}_2(\dsR_+)}+|r_n|_{\dot{\operatorname{H}}^{\theta+\frac{1}{2}-\sfH}_2(\dsR_+)}\right]^2<\infty
\end{equation}
implies Hölder continuity of $u(t)$ in time of order $\theta$.
Finally from \textbf{(e)} we obtain by interpolation
\[|e_n(\xi)-e_n(\eta)|\leq C|\xi-\eta|^\theta\mu_n^{\theta/2},\]
hence
\begin{equation}\label{sigma3}
  \sigma_3:=\sum_{n=1}^\infty\gamma_n\mu_n^\theta|r_n|^2_{\dot{\operatorname{H}}^{\frac{1}{2}-\sfH}_2(\dsR_+)}<\infty
\end{equation}
yields Hölder-continuity of $u(t,\xi)$ in space $\xi$ of order
$\theta$. Therefore the goal is to estimate the
$\dot{\operatorname{H}}^{\theta+\frac{1}{2}-\sfH}_{2}(\dsR_+)$-norms
of $r_n$, where the functions $r_n(t)$ are the fundamental solutions
of the scalar problems
\begin{equation}\label{scalar}
  r_n+\mu_n g_\alpha*r_n=g_\beta.
\end{equation}
This will be done by the following Lemma.

\begin{Lemma}\label{Lemma_B22}
  Suppose $\alpha\in(0,2)$, $\beta>0$, $\theta\in[0,1]$, and let
  $r_n(t)$ denote the solution of (\ref{scalar}). Then
  \begin{equation*}
    |r_n|^2_{\dot{\operatorname{H}}^{\theta+\frac{1}{2}-\sfH}_{2}(\dsR_+)}\leq
    C_{\alpha,\beta,\theta}\mu_n^{\frac{2(1-\beta+\theta-\sfH)}{\alpha}},\quad
    \mu_n>0,
  \end{equation*}
  whenever $\beta\in(1-\sfH+\theta,1-\sfH+\alpha)$.
\end{Lemma}

\noindent\textbf{Proof.} Again we extend the functions $r_n$
trivially on negative halfline. Let $\sfH\in(\tfrac{1}{2},1)$. We
first consider the case $\theta=0$. Then by the Paley-Wiener
theorem, $\calI^{\sfH-\frac{1}{2}}r_n\in L_2(\dsR)$ if and only if
$\hat{\calI^{\sfH-\frac{1}{2}}r_n}\in\mathscr{H}_2(\mathds{C}_+)$,
the Hardy space of exponent 2 and
$|\calI^{\sfH-\frac{1}{2}}r_n|_{L_2(\dsR)}=(1/\sqrt{2\pi})|\hat{\calI^{\sfH-\frac{1}{2}}r_n}|_{\mathscr{H}_2(\mathds{C}_+)}$.
Applying the Paley-Wiener theorem one more time, it suffices to show
that $\mathcal{F}(\calI^{\sfH-\frac{1}{2}}r_n)\in L_2(\dsR)$. Now we
may use (\ref{Fourier}) to compute
\[\int_\dsR \left|\mathcal{F}(\calI^{\sfH-\frac{1}{2}}r_n)(\rho)\right|^2\dd\rho\leq \int_0^\infty
\left[\frac{\rho^\alpha}{\rho^\beta(\rho^\alpha+\mu_n)}\right]^2\rho^{-2\sfH+1}\dd\rho\]
and we have seen in (\ref{rnnorm}) that the right integral converges
if and only if $\beta\in(1-\sfH,1-\sfH+\alpha)$. In case
$\theta\not=0$, observe that
$|\cdot|_{L_2(\dsR)}+|D^\theta\cdot|_{L_2(\dsR)}$ defines an
equivalent norm in $\operatorname{H}^\theta_{2}(\dsR)$, hence
replacing $\beta$ by $\beta-\theta$ the result follows by
Plancherel's Theorem. For $\sfH\in(0,\tfrac{1}{2})$ one may proceed
as above with replacing
$\calI^{\sfH-\frac{1}{2}}$ by $\calD^{\frac{1}{2}-\sfH}$.\qed\\

Now we are in the position to state our result on
(\ref{fracprob}).\\

\begin{Theorem}\label{Theorem2}
  Let $\alpha\in(0,2)$, $\beta>0$, $\theta\in[0,1]$ such that
  $\beta\in(1-\sfH+\theta,1-\sfH+\alpha)$. Assume that \textbf{(A)} and
  \textbf{(B)} are satisfied.
  \begin{enumerate}[(i)]
    \item If
    \[\sum_{n=1}^\infty
    \gamma_n\mu_n^{\frac{2(1-\beta-\sfH)}{\alpha}}<\infty,\] then the solution $u$
    of (\ref{fracprob}) exists and belongs to
    $C_b(\dsR_+;L_2(\Omega;\calH))$.
    \item If \[\sum_{n=1}^\infty
    \gamma_n\mu_n^{\frac{2(1-\beta+\theta-\sfH)}{\alpha}}<\infty,\] then $u\in
    C_b^\theta(\dsR_+;L_2(\Omega;\calH))$.
    \item If $\calH=L_2(G)$, \textbf{(e)} holds, and \[\sum_{n=1}^\infty
    \gamma_n\mu_n^{\frac{2(1-\beta-\sfH)+\alpha\theta}{\alpha}}<\infty,\] then
    $u\in C_b(\dsR_+;C^\theta(G;L_2(\Omega)))$.\\
  \end{enumerate}
\end{Theorem}

\noindent\textbf{Proof.} Use Lemma \ref{Lemma_B22} to estimate the
quantities $\sigma_i$, $i=1,2,3$, arising in (\ref{sigma1}),
(\ref{sigma2}) and (\ref{sigma3}), respectively. \qed

\noindent\begin{Rem}\
  \begin{enumerate}[(i)]
    \item In case $\beta=1$ and $\alpha>\sfH$, Theorem
    \ref{Theorem2} (i) coincides with the result of Bonaccorsi
    \cite[Theorem 4.8 (iii)]{Bon07}.
    \item Sufficient conditions for the existence of a mild solution
    of (\ref{fracprob}) in case $\beta=1$ and $\alpha\leq\sfH$ can
    be found in \cite[Theorem 4.8 (i),(ii)]{Bon07}.\\
  \end{enumerate}
\end{Rem}

\noindent\textbf{Example.} Let $\calH=L_2(0,\pi)$, $A=A_0^m$, where
$A_0=-(\dd/\dd x)^2$ with domain
$D(A_0)=\operatorname{H}^2_2(0,\pi)\cap\oH_2^1(0,\pi)$. It is
obvious that $A$ is subject to Hypothesis \textbf{(A)} and it is
well known that eigenvalues of $A$ are $\mu_k=k^{2m}$ for
$k\in\dsN$. The covariance $Q$ is given by its spectral
decomposition
\[Qx=\sum_{k=1}^\infty \gamma_k (x|e_k)e_k,\] with
$(\gamma_k)_{k\in\dsN}\subset(0,1]$ such that
$\sum_{k=1}^\infty\gamma_k<\infty$. For our example we choose
$\gamma_k=k^{-l}$, $l>1$, and we obtain
\[\begin{split}
  \sum_{k=1}^\infty
  \gamma_k\mu_k^{\frac{2(1-\beta-\sfH)}{\alpha}}<\infty\quad&\Longleftrightarrow\quad
  \beta>1-\sfH-\frac{\alpha(l-1)}{4m};\\
  \sum_{k=1}^\infty\gamma_k\mu_k^{\frac{2(1-\beta+\theta-\sfH)}{\alpha}}<\infty
  \quad&\Longleftrightarrow\quad
  \beta>1-\sfH+\theta-\frac{\alpha(l-1)}{4m};\\
  \sum_{k=1}^\infty
  \gamma_k\mu_k^{\frac{2(1-\beta-\sfH)+\alpha\theta}{\alpha}}<\infty
  \quad&\Longleftrightarrow\quad
  \beta>1-\sfH + \frac{\alpha\theta}{2} -\frac{\alpha(l-1)}{4m}.
\end{split}\]
Obviously the latter series converge for all
$\beta\in(1-\sfH+\theta,1-\sfH+\alpha)$, hence Theorem
\ref{Theorem2} applies independently from the choice of $l$ and $m$.
Observe that the spatial regularity is better than in time and that
for $\sfH\in(\tfrac{1}{2},1)$ the regularity in space and in time is
better than in case $\sfH=\tfrac{1}{2}$. On the other hand
regularity degrade for $\sfH\in(0,\tfrac{1}{2})$.\\

We conclude with a brief discussion of the case $\alpha=2$. Then
\begin{equation*}
  \hat{r}_n(\lambda)=\frac{\lambda^{2-\beta}}{\lambda^2+\mu_n},\qquad n\in\dsN,
\end{equation*}
hence there are poles $\pm i\sqrt{\mu_n}$ on the imaginary axis, and so Lemma \ref{Lemma_B22} is not valid in this case.
Therefore we proceed differently. It is shown in \cite{WNP}, that if $\tfrac{1}{2}<\beta<3$ one obtains with the aid of the complex inversion formula for the Laplace transform
\begin{equation*}
  r_n(t)=\mu_n^{\frac{1-\beta}{2}}\left[\sin\left(\sqrt{\mu_n}t+\frac{(2-\beta)\pi}{2}\right) - \frac{1}{\pi}\sin((2-\beta)\pi)\int_0^\infty e^{-\sqrt{\mu_n}t\tau}\frac{\tau^{2-\beta}\dd\tau}{1+\tau^2}\right],
\end{equation*}
where $t>0$. This formula shows in particular, that for every $t\in(0,T)$ it is $|r_n(t)|\leq c_T\mu_n^{\frac{1-\beta}{2}}$, where the constant $c_T$ (in the sequel generic) may depend on $T$.
%In case $\sfH<\tfrac{1}{2}$ this gives
%\begin{equation*}
%  |r_n|_{\dot{\operatorname{H}}^{\frac{1}{2}-\sfH}_{2}(0,T)}\leq |r_n|_{\operatorname{H}^{\frac{1}{2}-\sfH}_{2}(0,T)}\leq c|r_n|^{\frac{1}{2}+\sfH}_{L_2(0,T)}|r_n|^{\frac{1}{2}-\sfH}_{\operatorname{H}^1_2(0,T)}\sim c_T\mu_n^{\frac{3-2\beta}{4}-\frac{\sfH}{2}},
%\end{equation*}
%as $\mu_n\to\infty$ for any fixed $T>0$ and for all $n\in\dsN$.
Thanks to $L_{\frac{1}{\sfH}}\hookrightarrow \dot{\operatorname{H}}^{\frac{1}{2}-\sfH}_{2}$ we have in case $\sfH>\tfrac{1}{2}$
\begin{equation*}
  |r_n|_{\dot{\operatorname{H}}^{\frac{1}{2}-\sfH}_{2}(0,T)}\leq c_T |r_n|_{L_{\frac{1}{\sfH}}(0,T)}\leq c_T\mu_n^{\frac{1-\beta}{2}},
\end{equation*}
for any fixed $T>0$ and for all $n\in\dsN$.
The condition for local existence in the case $\alpha=2$ and $\sfH>\tfrac{1}{2}$ is now immediate and reads as
\begin{equation*}
  \sum_{n=1}^\infty\gamma_n\mu_n^{\frac{1-\beta}{2}}<\infty.
\end{equation*}
Note, that this is not the limiting case of Theorem \ref{Theorem2} (i) as $\alpha\to 2$.\\

\noindent\textbf{Acknowledgement:} We are grateful to J. Prüss and
W. Grecksch for some valuable suggestions which contributed to the
completion of this paper. Moreover we thank the anonymous referees for detailed reading and helpful comments.\\

\nocite{TTV03,Tud05,DJP06}
\bibliographystyle{plain}
\bibliography{fwn_pert}

\begin{thebibliography}{10}

\bibitem{BHOZ08}
F.~Biagini, Y.~Hu, B.~{\O}ksendal, and T.~Zhang.
\newblock {\em Stochastic calculus for fractional {B}rownian motion and
  applications}.
\newblock Probability and its Applications (New York). Springer-Verlag London
  Ltd., London, 2008.

\bibitem{Bon07}
S.~Bonaccorsi.
\newblock Volterra equations perturbed by a {G}aussian noise.
\newblock In {\em Progress in Probability}, volume~59, pages 37--55. Birkhäuser
  Verlag Basel/Switzerland, 2007.

\bibitem{CLE1}
Ph. Clément and G.~DaPrato.
\newblock {S}ome results on stochastic convolutions arising in volterra
  equations perturbed by noise.
\newblock {\em Atti Accad. Naz. Lincei (9) Mat. Appl.}, 7(3):147--153, 1996.

\bibitem{CLE2}
Ph. Clément and G.~DaPrato.
\newblock {W}hite noise perturbations of the heat equation in materials with
  memory.
\newblock {\em Dynam. Systems Appl.}, 6(4):441--460, 1997.

\bibitem{WNP}
Ph. Clément, G.~DaPrato, and J.~Prüss.
\newblock White noise perturbation of the linear parabolic viscoelasticity.
\newblock {\em Rend. Istit. Mat. Univ. Trieste}, 29(1-2):207--220, 1997.

\bibitem{StochEqua}
G.~DaPrato and J.~Zabczyk.
\newblock {\em {S}tochastic equations in infinite dimensions}.
\newblock Encyclopedia of Mathematics and its Applications. Cambridge
  University Press, 1992.

\bibitem{DJP06}
T.~E. Duncan, J.~Jakubowski, and B.~Pasik-Duncan.
\newblock Stochastic integration for fractional {B}rownian motion in a
  {H}ilbert space.
\newblock {\em Stoch. Dyn.}, 6(1):53--75, 2006.

\bibitem{KaSh91}
I.~Karatzas and S.~E. Shreve.
\newblock {\em Brownian motion and stochastic calculus}, volume 113 of {\em
  Graduate Texts in Mathematics}.
\newblock Springer-Verlag, New York, second edition, 1991.

\bibitem{KwWo92}
S.~Kwapie{\'n} and W.~A. Woyczy{\'n}ski.
\newblock {\em Random series and stochastic integrals: single and multiple}.
\newblock Probability and its Applications. Birkh\"auser Boston Inc., Boston,
  MA, 1992.

\bibitem{MaNe68}
B.~B. Mandelbrot and J.~W. Van~Ness.
\newblock Fractional {B}rownian motions, fractional noises and applications.
\newblock {\em SIAM Rev.}, 10:422--437, 1968.

\bibitem{Mon}
S.~Monniaux and J.~Prüss.
\newblock {A} theorem of the {D}ore-{V}enni type for noncommuting operators.
\newblock {\em Trans. Amer. Math. Soc.}, 349(12):4787--4814, 1997.

\bibitem{Nua06}
David Nualart.
\newblock {\em The {M}alliavin calculus and related topics}.
\newblock Probability and its Applications (New York). Springer-Verlag, Berlin,
  second edition, 2006.

\bibitem{FBM-Int}
V.~Pipiras and M.~S. Taqqu.
\newblock {I}ntegration queations related to fractional {B}rownian motion.
\newblock {\em Probab. Theory Relat. Fields}, 118(2):251--291, 2000.

\bibitem{EvolInt}
J.~Prüss.
\newblock {\em {E}volutionary integral equations and applications}, volume~87
  of {\em Monographs in Mathematics}.
\newblock Birkhäuser, 1993.

\bibitem{Samko}
S.~G. Samko, A.~A. Kilbas, and O.~I. Marichev.
\newblock {\em Fractional integrals and derivatives}.
\newblock Gordon and Breach Science Publishers, 1994.

\bibitem{TTV03}
S.~Tindel, C.~A. Tudor, and F.~Viens.
\newblock Stochastic evolution equations with fractional {B}rownian motion.
\newblock {\em Probab. Theory Related Fields}, 127(2):186--204, 2003.

\bibitem{Tri83}
H.~Triebel.
\newblock {\em Theory of function spaces}, volume~78 of {\em Monographs in
  Mathematics}.
\newblock Birkh\"auser Verlag, Basel, 1983.

\bibitem{Triebel}
H.~Triebel.
\newblock {\em {I}nterpolation theory, function spaces, differential
  operators}.
\newblock Johann Ambrosius Barth, Heidelberg, 1995.

\bibitem{Tud05}
C.~A. Tudor.
\newblock It\^o formula for the infinite-dimensional fractional {B}rownian
  motion.
\newblock {\em J. Math. Kyoto Univ.}, 45(3):531--546, 2005.

\end{thebibliography}
\end{document}